\documentclass[review]{elsarticle}

\usepackage{lineno,hyperref}
\modulolinenumbers[5]

\journal{}

\usepackage{graphicx}
\usepackage[mathscr]{eucal}
\usepackage{amsmath}
\usepackage{amssymb}
\usepackage{graphics}

\newtheorem{theorem}{Theorem}[section]

\newtheorem{definition}[theorem]{Definition}

\begin{document}

\begin{frontmatter}

\title{Obstacle avoiding patterns and \\cohesiveness of fish school}%\tnoteref{mytitlenote}}
%\tnotetext[mytitlenote]{Fully documented templates are available in the elsarticle package on \href{http://www.ctan.org/tex-archive/macros/latex/contrib/elsarticle}{CTAN}.}

%% Group authors per affiliation:
%\author{Linh Thi Hoai Nguyen\fnref{myfootnote1}, Ta Viet Ton}
%\address{Department of Information and Physical Sciences, Graduate School of Information Science and Technology, Osaka University, Suita, Osaka 565-0871, Japan}
%\fntext[myfootnote]{The work of the second author is supported by JSPS KAKENHI Grant Number 20140047.}

%% or include affiliations in footnotes:
\author[address1]{Linh Thi Hoai Nguyen\corref{mycorrespondingauthor}}
\cortext[mycorrespondingauthor]{Corresponding author}
\ead{nth.linh@ist.osaka-u.ac.jp}

\author[address2]{Vi$\hat{\d{e}}$t T$\hat{\rm o}$n T\d{a}}
%\fntext[fund]{The work of the second author is supported by JSPS KAKENHI Grant Number 20140047.}
\ead{taviet.ton@ist.osaka-u.ac.jp}

\author[address3]{Atsushi Yagi}
\ead{atsushi-yagi@ist.osaka-u.ac.jp}

\address[address1]{Department of Information and Physical Sciences, Graduate School of Information Science and Technology, Osaka University, Suita, Osaka 565-0871, Japan}
\address[address2]{Promotive Center for International Education and Research of Agriculture\\  
Faculty of Agriculture,  
 Kyushu University, Nishi-ku, Fukuoka 812-8581, Japan}
\address[address3]{Department of Applied Physics, Graduate School of Engineering, Osaka University, Suita, Osaka 565-0871, Japan}

\begin{abstract}
This paper is devoted to studying obstacle avoiding patterns and cohesiveness of fish school. First, 
we introduce a model of stochastic  differential equations (SDEs) for describing the process of   fish school's obstacle avoidance.
Second, on the basis of the model we find obstacle avoiding patterns. Our observations show that there are clear four  obstacle avoiding patterns, namely, Rebound, Pullback, Pass and Reunion, and Separation. Furthermore, the emerging patterns change when parameters change. Finally, we present a scientific definition for fish school's  cohesiveness that will be an internal property characterizing the strength of fish schooling. There are then evidences that 
 the school cohesiveness can be measured through obstacle avoiding patterns.
\end{abstract}


\end{frontmatter}

%\linenumbers

\section{Introduction}\label{intro}

Fish schooling, one of animal swarming, is a commonly observed phenomenon that is coherently performed by integration of interactions among constituent fish. This remarkable phenomenon has already attracted interests of researchers from diverse fields including biology, physics, mathematics, computer engineering (see \cite{Aoki1982,Bonabeau1999,Camazine2001,Gunji1999,Huth1992,Nguyen2014,Olfati2006,Olfati2003,Reynolds1987,Ta,Uchitane2012}).

%For example, biologists empirically report how a number of small fish behave like a single living thing maintaining their school without any collision with each other. Physicists present their insight on the mechanism how such a coherent behavior is created by interactions among individual fish and write down these interactions by mathematical formulas. Mathematicians then deduce interesting information from the equation to clarify the properties that the fish school has inherently. Such an equation can also give useful suggestions to computer engineers to design various computer systems which run automatically and spontaneously.

Let us recall here some researches in the literature therein. In 2001, Camazine et al.  \cite[Chapter 11]{Camazine2001}  presented an idea on the basis of experimental results (Aoki \cite{Aoki1982}, Huth-Wissel \cite{Huth1992}, and Warburton-Lazarus \cite{WL}) that individual fish may act following the behavioral rules:
\begin{enumerate}
\item [(a)] The school has no leaders and each fish follows the same behavioral rules.
\item [(b)] To decide where to move, each fish uses some form of weighted average of the position and orientation of its nearest neighbors.
\item [(c)] There is a degree of uncertainty in the individual's behavior that reflects both the imperfect information-gathering ability of a fish and the imperfect execution of the fish's actions.
\end{enumerate}
Their insight is that these local rules can altogether create the coherent behavior of fish school.

 Vicsek et al. \cite{Vicsek1995} modeled the movement as self-driven particles obeying some  difference equations. They assumed that each individual is driven with a sum of an absolute velocity  and an averaged velocity of nearby particles together with some random perturbations. Oboshi et al. \cite{Oboshi2002} also modeled  schooling by some difference equations setting a rule that each fish choose one way of action among four possibilities according to a distance to the closest mate. Meanwhile, Olfati-Saber \cite{Olfati2006} and D'Orsogna et al. \cite{DOrsogna2006} independently presented differential equation models, but deterministic ones, utilizing the generalized Morse function and attractive/repulsive potential functions, respectively. In \cite{Gunji1999}, Gunji et al.  considered dual interaction which produces territorial and schooling behavior. 
In \cite{Reynolds1987},  Reynolds introduced some simple behavioral rules of animals, which are similar to (a), (b), (c) but deterministic ones, and are schematic rather than physical. 
%In \cite{Reynolds1987},  Reynolds introduced some simple behavioral rules of animals which is schematic (rather than physical) and reconstructed their coherent behavior by computer simulations.
%(Similar assumptions, but deterministic ones, were also introduced by Reynolds \cite{Reynolds1987}.) 

A stochastic differential equation (SDE) model describing the process of schooling was presented in  \cite{Uchitane2012}, where we used the behavioral rules  (a), (b), and (c)  above. We then utilized the model for developing  quantitative arguments on fish schooling in \cite{Nguyen2014}.

%we quote Chang et al. \cite{Chang2003},  Hettiarachchi and  Spears \cite{Hettiarachchi2005}, Olfati-Saber and Murray \cite{Olfati2003}, and   Reynolds \cite{Reynolds1987}.

In the real world, the environment surrounding fish  school  often includes other components such as obstacles, food resources, predators, etc. In those situations, fish exhibit more complex, parallel movements such as obstacle avoidance, food finding, escaping from predator. It is evident that when a school of fish is tackled by obstacles or is hunted by predators, fish individually react quickly for avoiding obstacles or predators.  

%To authors' knowledge, there are so far many papers handling  the swarming behavior in the presence of  predators (see \cite{Axelsen2001,Oboshi2002,Partridge1982}, to name a few) but only a few papers dealing with it  in the presence of obstacles \cite{Chang2003,Hettiarachchi2005,Olfati2003,Reynolds1987}. 

 Olfati-Saber and Murray \cite{Olfati2003} developed the method of Reynolds  \cite{Reynolds1987} by introducing a dynamic graph of agents in the presence of multiple obstacles. The agents are split into several groups while approaching the obstacles. After passing all obstacles, they rejoin into a single group.
Chang et al. \cite{Chang2003}  introduced techniques of using gyroscopic forces for multi-agent systems by which the agents perform collision avoidance toward obstacles. A similar result  has been shown (i.e., agents are separated into some clusters and then rejoin into a single flock). In the meantime,  Hettiarachchi and  Spears  \cite{Hettiarachchi2005} used virtual physical forces in composing a swarm system of robots moving toward a goal through  obstacle fields. Robots may collide with obstacles but then they can still move toward the goal.

In the meantime, a concept concerning  animal swarming  or grouping, namely cohesiveness has already been introduced since 1930s. The study during long years seems to show that it is not an easy problem to define a concept of the cohesiveness precisely and consistently. It has been conceptualized in various ways, but each was  based on  intuitive assumptions and interpretations. 

For instance,  Moreno and Jennings \cite{Moreno1937} defined cohesiveness as the forces holding the individuals within the group to which they belong. French \cite{French1941} noted that the group exists as a balance between cohesion and disruptive forces. Not until 1950 was a systematic theory of group cohesiveness constructed by Festinger et al. \cite{Festinger1950}. Their definition of cohesiveness is ``We shall call the total field of forces which act on members to remain in the group the ``cohesiveness'' of that group''.  
Gross and Martin \cite{Gross1952} claimed that this definition is inadequate, and they proposed an alternative definition as the resistance of group to disruptive forces. Contemporary works almost characterize group cohesion  in the same way (see \cite{Hogg1992}). Carron \cite{Carron1980} defined cohesiveness to be the adhesive property of group. Schachter et al. \cite{Schachter1951} found that interpersonal attraction is the cement binding group members together. For the general relationship between cohesiveness and group performance, we refer the reader to \cite{Beal2003, Laurel1988, Mullen1994}.

We can however find a point of view which is common in those definitions. That is the bond linking group members to  others and to the group as a whole.  We believe that this common point of view may be a key feature for all intercommunicated multi-agent systems.

The objective of the present paper is two-folds: namely, studying the fish schooling from a viewpoint of pattern formation of biological systems and introducing a scientific definition of the  school cohesiveness. 

For the first objective, obstacle avoiding patterns of fish school are studied by newly introducing a behavioral rule for avoidance and adding its effect to our model in \cite{Uchitane2012}. It is then observed that there are at least four obstacle avoiding patterns of school, i.e.,  Rebound, Pullback, Pass and Reunion, and Separation which are performed  just by tuning modeling parameters. 

For the second objective, we consider school cohesiveness as its ability  to form and maintain the school structure against the white noises affecting the school. It is therefore defined as an internal nature of the school, independent of external effects. We then show how internal parameters contribute to the school's cohesiveness. Furthermore, our results  suggest  a very interesting correlation between the degree of  cohesiveness and the four  avoidance patterns. 

The outline of this paper is as follows. Section \ref{MathematicalModels} gives model description. We first recall the SDE model for fish schooling introduced in \cite{Uchitane2012}, then newly inoculate a mechanism for obstacle avoiding into it. Section \ref{AvoidingObstaclePattern} presents four obstacle avoiding patterns. We thereafter investigate how these patterns change as the  modeling parameters are tuned. Section \ref{Cohesiveness} explores fish school cohesiveness. A scientific definition   and  measurement of  cohesiveness are introduced. The relationship between  avoidance patterns and school cohesiveness is then investigated.  The paper concludes with some discussions of Section \ref{Conclusions}.

\section{Model description}\label{MathematicalModels}

In \cite{Uchitane2012}, we  introduced a SDE model of the form
\begin{equation}\label{eq1}
\begin{cases}
d{\mathbf x}_i(t)&={\mathbf v}_idt+\sigma_idw_i(t),\\
d{\mathbf v}_i(t)&=\Big\{-\alpha\sum\limits_{j=1,j\ne i}^N\left(\dfrac{r^p}{\|{\mathbf x}_i-{\mathbf x}_j\|^p}-\dfrac{r^q}{\|{\mathbf x}_i-{\mathbf x}_j\|^q} \right)({\mathbf x}_i-{\mathbf x}_j)\\
  &\qquad -\beta\sum\limits_{j=1,j\ne i}^N\left(\dfrac{r^p}{\|{\mathbf x}_i-{\mathbf x}_j\|^p}+\dfrac{r^q}{\|{\mathbf x}_i-{\mathbf x}_j\|^q} \right)({\mathbf v}_i-{\mathbf v}_j)\\
&\qquad +F_i({\mathbf x}_i,{\mathbf v}_i)\Big\}dt,\hfill i=1,2,\ldots,N,
\end{cases}
\end{equation}
for an $N$-fish system moving in the space $\mathbb R^d$ ($d=2,3$). Here, ${\mathbf x}_i(t)$ and ${\mathbf v}_i(t)$ denote position and velocity  of the $i$-th fish at time $t$, respectively. And  $\|\cdot\|$ denote the Euclidean norm of a vector, hence $\|{\mathbf x}_i-{\mathbf x}_j\|$ represents the distance between the $i$-th and the $j$-th fish.

We regarded each fish as a moving particle in $\mathbb R^d$. The first equation is a stochastic equation for the unknown ${\mathbf x}_i(t)$, where $\sigma_idw_i$  denotes noise resulting from the imperfectness of information-gathering and action of the $i$-th fish. In fact, $w_i(\cdot) (i=1,2,\dots,N)$ are independent $d$-dimensional  Brownian motions on some probability space. %defined in a probability space $(\Omega,\mathcal F,{\mathbb P})$ with filtration $\{\mathcal F_t\}_{0\leqslant t<\infty}.$  

The second equation is a deterministic equation for ${\mathbf v}_i(t)$, where $1<p<q<\infty$ are fixed exponents; $\alpha$ and $\beta$ are positive coefficients of attraction and velocity matching among fish, respectively. And $r>0$ is a fixed number. 
If $\|{\mathbf x}_i-{\mathbf x}_j\|>r$ then the $i$-th fish moves toward the $j$-th.  To the contrary, if  $\|{\mathbf x}_i-{\mathbf x}_j\|<r$ then the $i$-th fish acts in order to avoid collision with the $j$-th, $r$ being thereby a critical distance (for details, see \cite{Uchitane2012}).

Velocity matching of the $i$-th fish to the $j$-th also has a similar weight depending on the distance $\|{\mathbf x}_i-{\mathbf x}_j\|$. Degree of matching is higher when $\|{\mathbf x}_i-{\mathbf x}_j\|<r$ for urgent reaction to avoid collision. When $p$ is large,  the attraction range is short. Meanwhile, when  $p$ is small,  fish have ability to attract each other against a long distance.  The exponent $p$ thereby relates to a degree of how far the attraction reaches. Finally, the function $F_i({\mathbf x}_i,{\mathbf v}_i)$ stands for an external force function acting on the $i$-th fish.  

The system  \eqref{eq1} was studied analytically and numerically in the subsequent paper \cite{Nguyen2014} of  \cite{Uchitane2012}.   Quantitative arguments on fish schooling were also developed.
Main advantages of using SDE models like \eqref{eq1} over others may be the convenience of mathematical technique. One can utilize the well developed theory of SDEs including  numerical methods (see \cite{Kloeden2005}). Its flexibility may be another advantage. As seen below, it is very easy to modify the model in the existence of an obstacle simply by introducing a suitable external force functions $F_i(\mathbf x_i,\, \mathbf v_i)$ in \eqref{eq1}. 

%Similar results might be, however, obtained by using other models too.

We now want to  study \eqref{eq1} from the viewpoint of obstacle avoiding. 
For this purpose, let us put a global obstacle with central point ${\mathbf x}_C$ and radius $\rho>0$ in $\mathbb R^d$ and assume that fish move in the domain $\Omega=\{{\mathbf x}\in\mathbb R^d:\|{\mathbf x}-{\mathbf x}_C\|>\rho\}$. The surface of the obstacle is denoted by $S=\{\mathbf x\in\mathbb R^d:\|{\mathbf x}-{\mathbf x}_C\|=\rho\}$. In addition to the rules (a)-(c) in Introduction, we introduce a behavioral rule in order to avoid collision with the obstacle:
\begin{itemize}
\item[(d)] Each fish executes an action for avoiding the obstacle according to a reflection law of velocity with distance depending weight.
\end{itemize}

On the basis of this rule, we are able to  formulate the external force function $F_i(x_i,v_i)$ which influences the motion of the $i$-th individual. For a given $({\mathbf x}_i,{\mathbf v}_i)$, where $\|{\mathbf x}_i-{\mathbf x}_C\|>\rho$, let $l_i$ be a ray with origin ${\mathbf x}_i$ and direction ${\mathbf v}_i$, i.e.,   
$$l_i=\{{\mathbf x}\in\mathbb R^d:{\mathbf x}={\mathbf x}_i+s{\mathbf v}_i,\quad 0\leqslant s<\infty\}.$$
When $l_i$  intersects $S$, we define ${\rm Rf}({\mathbf x}_i,{\mathbf v}_i)$ as the reflection vector ${\mathbf u}_i$ of ${\mathbf v}_i$ at ${\mathbf y}_i$ with respect to the tangential plane of $S$ at ${\mathbf y}_i$, where ${\mathbf y}_i$ is the first intersection of $l_i$ to $S$. In fact, ${\mathbf u}_i$ is a vector whose opposite equals to the symmetric vector of ${\mathbf v}_i$ over the line connecting ${\mathbf y}_i$ and ${\mathbf x}_C$  on the plane generated by ${\mathbf v}_i$ and  ${\mathbf x}_C$.
When $l_i$ does not intersect $S$ (including the special case when ${\mathbf v}_i={\mathbf 0}$),  we put  ${\rm Rf}({\mathbf x}_i,{\mathbf v}_i)={\mathbf v}_i$. Fig. 1 illustrates $({\mathbf x}_i,{\mathbf v}_i)$ and ${\rm Rf}({\mathbf x}_i,{\mathbf v}_i)$ in two-dimensional space. 
\begin{center}
\begin{figure}[h!]
	\includegraphics[width=0.8\textwidth]{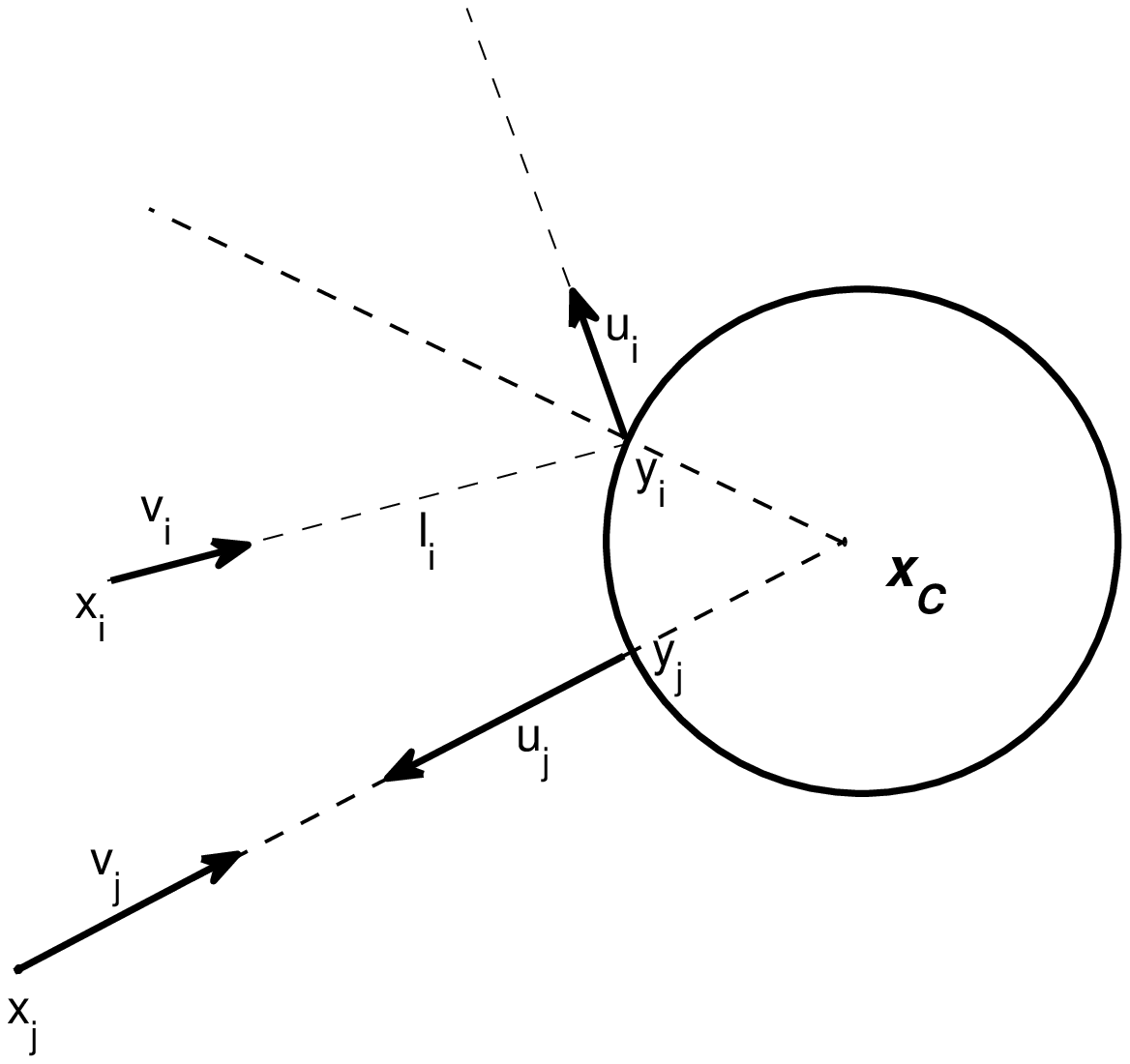}
	\caption{Reflection vector ${\rm Rf}(\bold x_i,\bold v_i)=\bold u_i$}
	\label{fig:1}   
\end{figure}
\end{center}
Analogously to the velocity matching, we  formulate 
\begin{align}\label{eq2}
F_i({\mathbf x}_i,{\mathbf v}_i)=-\gamma\left( \frac{R^P}{\|{\mathbf x}_i-{\mathbf y}_i\|^P}+\frac{R^Q}{\|{\mathbf x}_i-{\mathbf y}_i\|^Q}\right)[{\mathbf v}_i-{\rm Rf}({\mathbf x}_i,{\mathbf v}_i)],
\end{align}
where $1<P<Q<\infty$ are exponents, $R$ is a fixed distance, and $\gamma>0$ is a constant.

If $\|{\mathbf x}_i-{\mathbf y}_i\|<R$, then the $i$-th fish promptly reacts for matching its velocity to the reflection vector ${\rm Rf}({\mathbf x}_i,{\mathbf v}_i)$ to avoid collision with the obstacle. Meanwhile, if  $\|{\mathbf x}_i-{\mathbf y}_i\|>R,$ then the reaction to avoid the obstacle is less strong. If the ray $l_i$ does not meet $S$, the individual takes no reaction with the obstacle.

The model for fish schooling in obstacle domain is therefore of the form \eqref{eq1} with the external forces $F_i({\mathbf x}_i,{\mathbf v}_i)$  $(i=1,2,\ldots,N)$ given by  \eqref{eq2}. 

%%%%%%%%%%%%%%%%%%%%%%%%%%%%%%%%%%%%%%%%%%%%%%%%%%%%%%%%%%%%%%%%%%%%%%%%%%%%%%%%
%%%%%%%%%%%%%%%%%%%%Cohesiveness and obstacle avoidance patterns Section     	%%%%%%%%%%%%%%%%%%%%%%%%%%%%%%%%%%%%%%%%%%%%%%%%%%%%%%%%%%%
\section{Obstacle avoiding patterns}\label{AvoidingObstaclePattern}
In this section, we  observe  four  avoidance patterns of fish school based on our model \eqref{eq1}, where the external force functions $F_i({\mathbf x}_i,{\mathbf v}_i)\, (i=1,\dots,N)$ are defined by  \eqref{eq2}. The effect of  control  parameters to these patterns is also investigated.

\subsection{$\varepsilon$-Graph and $\varepsilon,\theta$-Schooling}
Let us review  notions of $\varepsilon$-Graph and $\varepsilon,\theta$-Schooling that were introduced in the paper \cite{Nguyen2014}.
\begin{definition}[$\varepsilon$-Graph]  
Let ${\mathbf x}_i,{\mathbf v}_i$ $(i=1,2,\ldots,N)$ denote a solution to SDEs \eqref{eq1}. Let $\varepsilon>0$ be a fixed length. Regard the positions ${\mathbf x}_i(t)$ of  fish at each time $t$ as the vertices of a graph. Two vertices ${\mathbf x}_i(t)$ and ${\mathbf x}_j(t)$ are said to be connected by an edge of graph if and only if $\|{\mathbf x}_i(t)-{\mathbf x}_j(t)\|\leqslant \varepsilon$. Such a graph is called the $\varepsilon$-graph of group at time $t$ and is denoted by ${\rm G}_{\varepsilon}(t)$.
\end{definition}

Denote by $n_{\varepsilon}(t)$  the number of connected components of ${\rm G}_{\varepsilon}(t)$. When $n_{\varepsilon}(t)=1$,  the individuals  form a single group. Meanwhile, when  $n_{\varepsilon}(t)\geqslant 2$, they form $n_{\varepsilon}(t)$  sub-groups.

Denote by $\sigma {\rm V}(t)$ the variation of velocities:
$$\sigma{\rm V}(t)=\sqrt{\frac1N\sum\limits_{i=1}^N\|{\mathbf v}_i(t)-\bar {\mathbf v}(t)\|^2},\hspace{2cm} 0<t<\infty,$$
where $\bar {\mathbf v}(t)=\frac{1}{N}\sum\limits_{i=1}^N{\mathbf v}_i(t)$ is the average of all velocities of fish at time $t$.
\begin{definition}[$\varepsilon,\theta$-Schooling]\label{SchoolingDef}
 Let $\varepsilon>0$ and $\theta>0$ be given. If a solution $({\mathbf x}_i,{\mathbf v}_i)$ $(i=1,2,\ldots,N)$ to \eqref{eq1} satisfies $n_{\varepsilon}(t)=1$ and $\sigma{\rm V}(t)\leqslant \theta$ for all sufficiently large $t$, say all $t\geqslant T$ with some fixed time $T>0$, then the state of fish  is said to be in $\varepsilon,\theta$-schooling.
\end{definition}

When the distance $\varepsilon>0$ and the tolerance $\theta>0$  are specified,  we simply say ``in schooling" instead of ``in $\varepsilon,\theta$-schooling".

%------------------------------- Avoidance patterns --------------------------------------------------------------
\subsection{Obstacle avoiding patterns}    \label{Patterns}
Set  parameters as  $d=2$, $N=20$, $\alpha=\beta=\gamma=1$. The exponent $p$ is tuned from 2 to 4 keeping always the relation $q=p+1$ and $\sigma_i=0$ for all $i$. The critical distance $r$ is set by $r=0.5$, and the radius of the obstacle is $\rho=1.2$. In addition, $\varepsilon=0.5$ and $\theta=10^{-6}$.

%Modeling parameters are chosen as $\varepsilon=r=0.5$ and $\theta=10^{-6}$. 

By performing preliminary computations, we first set a stationary state which is in $\varepsilon,\theta$-schooling. This state is set as the initial position of our computations. The distance from the center of the school to the center of the obstacle is 3.5, and the line connecting these two centers coincides with the horizontal axis. The initial velocities are ${\mathbf v}_i(0)=(1.75,\, 0)$ for every $i$. Parameters for obstacle avoidance are set as $P=p$, $Q=q$, $R=r$. The school is thereby oriented toward the obstacle and will  strike on it after a while.

\begin{figure*}[h!]
	\includegraphics[width=1\textwidth]{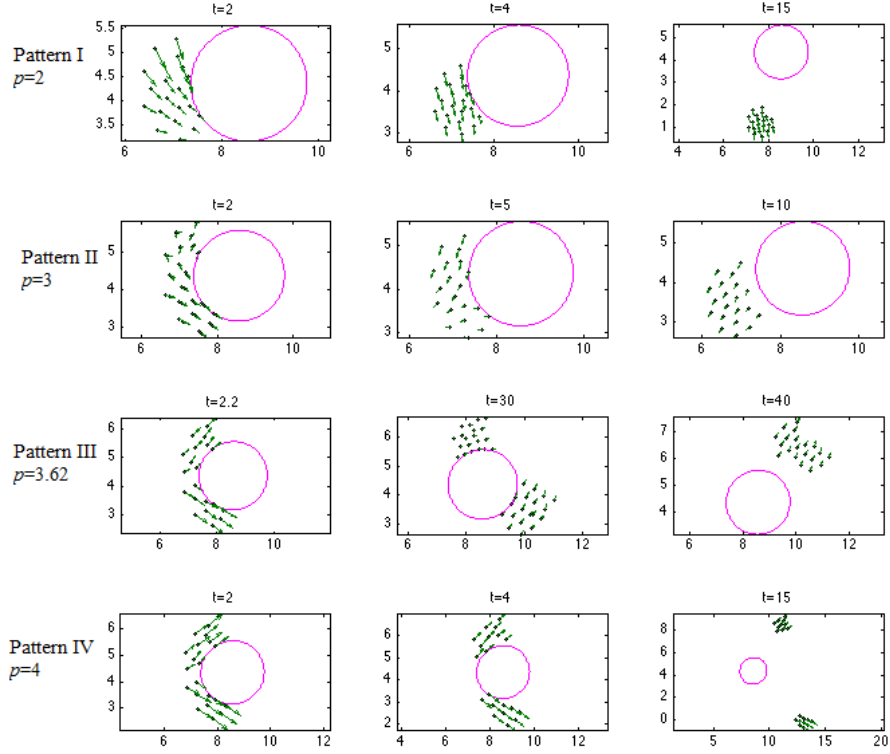}
	\caption{Obstacle avoiding patterns in two-dimensional space: Rebound,  Pullback,  Pass and Reunion, and  Separation}
	\label{fig:5} 
\end{figure*}

Fig.~\ref{fig:5}  shows the numerical results for $p=2,\,3,\,3.62,\, 4$. Four different kinds of  avoiding patterns  are found. We  will call them,  Rebound,  Pullback,  Pass and Reunion, and  Separation, respectively.

Let us describe these four  patterns of fish schooling.

{\bf Pattern I (Rebound)}: The fish keep schooling throughout the obstacle avoiding process and the school rebounds off the obstacle. In order to keep being in schooling, they change their directions after the school touch the obstacle.

{\bf Pattern II (Pullback)}: The individuals are once separated while approaching the obstacle and stay around the surface of obstacle for a while. They then pull back off the obstacle to reform a school structure.

{\bf Pattern III (Pass and Reunion)}: The fish pass the obstacle by spliting to move along the obstacle surface.  After passing it they reunite into a single school.

{\bf Pattern IV (Separation)}: It is similar to Pattern III. But, after passing the obstacle,  the subgroups have their own direction.

Similar results are observed in the three-dimensional case, i.e.,   $d=3$ too. 
Numerical solutions to \eqref{eq1} show that there are similarly  four   avoidance patterns (I)-(IV).

\begin{figure*}[h!]
	\includegraphics[width=1\textwidth]{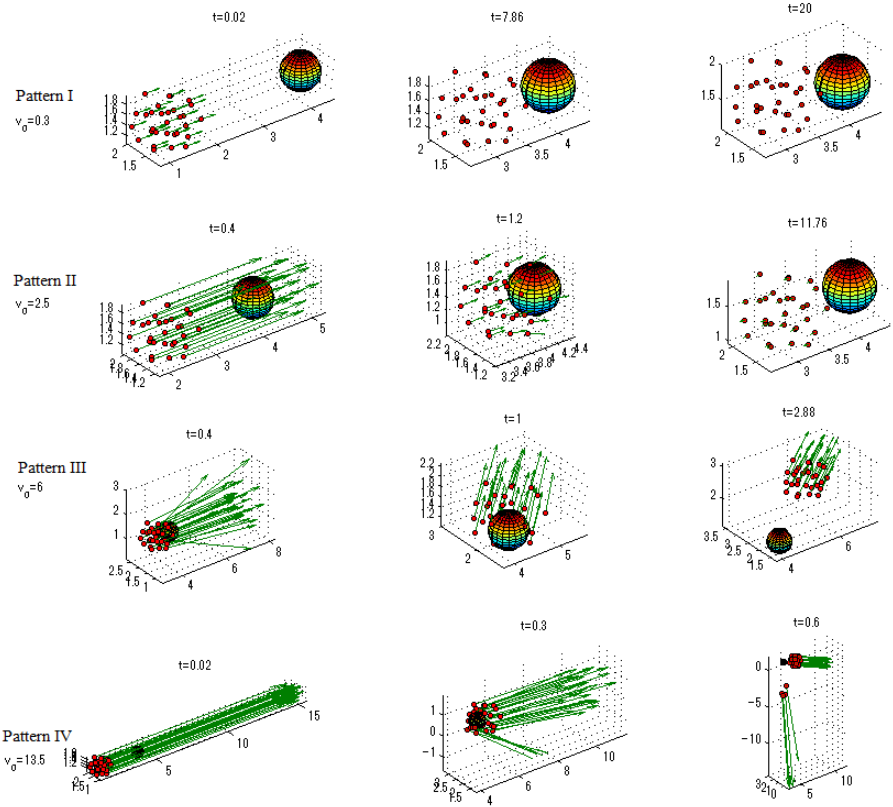}
	\caption{Obstacle avoiding patterns in three-dimensional space: Rebound,  Pullback,  Pass and Reunion, and  Separation}
	\label{3DPatterns} 
\end{figure*}

Fig. \ref{3DPatterns} represents four behavioral patterns of fish school while avoiding a static sphere obstacle in three-dimensional space. Four rows in this figure illustrate the time evolution of school starting from  different  initial speeds $v_0=0.3,\, 2.5, \, 6, 13.5$, respectively, while the other parameters are kept being constant in such a way that  $p=2$, $\alpha=\beta=\gamma=2$, $r=1$. Initial positions of fish are set   in $\varepsilon, \theta$-schooling, where $\varepsilon=1$ and $\theta=10^{-6}$.

%------------------------------- Exponents $p,q$  and patterns --------------------------------------------------------------
\subsection{Exponents $p,q$  and obstacle avoiding patterns}   \label{Exponents-patterns}
Set the parameters:  $d=2$, $N=20$, $\alpha=\beta=\gamma=1$, $r=0.5$, $\rho=1.2$,  $\varepsilon=0.5$ and $\theta=10^{-6}$ as Subsection  \ref{Patterns}. The initial data (velocity and position) are  also the same. 

The exponent $p$ is however finely tuned  from 1.001 to 8 with increment $\Delta p=0.001$ and the relation $q=p+1$ is always kept. In addition, $P=p$, $Q=q$, $R=r$.

Our numerical results are presented in Table 1. This table  highlights that there are critical values of $p$ at which the type  of  avoidance patterns changes from I to II, from II to III, and from III to IV, respectively. 

\begin{table}[h!]
\caption{Exponent $p$ and patterns}
\label{tab:1}       % Give a unique label
\begin{tabular}{lllll}
\hline\noalign{\smallskip}
$p$ & [1.001,\, 2.100] & [2.101,\, 3.371]&[3.372,\,3.497]&[3.498,\,8.000]  \\
\noalign{\smallskip}\hline\noalign{\smallskip}
Pattern & I & II&III&IV \\
\noalign{\smallskip}\hline
\end{tabular}
\end{table}

\subsection{Speeds of school  and obstacle avoiding patterns}
We set a fish group being in schooling. We investigate relations between the speed of school $\|\bar v\|$ and the type of performed pattern.

Table~\ref{tab:3} shows our numerical results in two-dimensional case. Here $N=20$, $\alpha=\beta=\gamma=1$, $p=2$, $r=0.5$, and the distance from the center of fish in schooling to the center of the obstacle is 3.5. The increment $\Delta \|v_0\|$ is also 0.001.

\begin{table}[h!]
\caption{Speeds of school and patterns}
\label{tab:3}       % Give a unique label
\begin{tabular}{lllll}
\hline\noalign{\smallskip}
$\|v_0\|$ & [0.001,\, 1.199] & [1.200,\, 2.589]&[2.590,\,4.866]&[4.867,\,20]  \\
\noalign{\smallskip}\hline\noalign{\smallskip}
Pattern & I & II&III&IV \\
\noalign{\smallskip}\hline
\end{tabular}
\end{table}

Of course, if the initial speed is too large, then collision may happen because the large velocity makes fish  have no enough time to adjust to avoid collision with other fish or with the obstacle. We do not consider this case. It is, however, interesting to know that  solutions to \eqref{eq1}, where $F_i (i=1,2,\dots,N)$ are given by \eqref{eq2}, can blow up in a finite time.

%In this case the fragmentation happens because the distance between some fish approach zeros.
%------------------------------- Critical distance $r$ and patterns --------------------------------------------------------------
\subsection{Critical distance $r$ and obstacle avoiding patterns}
Let tune the critical distance $r$ from 0.2 to 2.8 with increment $\Delta r=0.1$. Other parameters are $d=2$, $N=20$, $\alpha=\beta=\gamma=1$, $\rho=1.2$, $P=p=3$, $Q=q=4$, $\varepsilon=R=r$, $\theta=10^{-6}$. 

In order to set the initial positions, we perform preliminary computations for each $r$ in the free space $\mathbb R^2$, where  $F_i=-5{\mathbf v}_i$. These computations provide stationary states in $\varepsilon,\theta$-schooling for each $r$.

As pointed out by \cite{Nguyen2014}, the geometrical diameter $\delta$ 
 $$\delta=\max_{1\leqslant i\leqslant N}\|{\mathbf x}_i-\bar {\mathbf x}\|  \hspace{1cm} \text{with } \bar {\mathbf x}=(\sum_{i=1}^N {\mathbf x}_i)/N$$
 of $\varepsilon,\theta$-schooling depends on  $r$. It is thus natural to choose different radius of obstacle depending on $r$. 
In the present computations,  the radius of obstacle is set as $\rho=2\delta(r)/3$, where $\delta(r)$ is the school diameter corresponding to $r$. 
The distance from the center of school to the center of obstacle is 8, and the line connecting these two centers coincides with  the horizontal axis. The initial velocities are ${\mathbf v}_i(0)=(4,0)$ for all $i$.

Our numerical results  given in Table~\ref{tab:2} show that as $r$ increases, the type of patterns changes from IV to III, from III to II, and from II to I.

\begin{table}[h!]
\caption{Critical distance  and patterns}
\label{tab:2}       % Give a unique label
\begin{tabular}{lllll}
\hline\noalign{\smallskip}
$r$ & [0.2,\, 0.3] & [0.4,\, 0.5]&[0.6,\,2.0]&[2.1,\,2.8]  \\
\noalign{\smallskip}\hline\noalign{\smallskip}
Pattern & IV & III&II&I \\
\noalign{\smallskip}\hline
\end{tabular}
\end{table}

%%%%%%%%%%%%%%%%%%%%%%%%%%%%%%%%%%%%%%%%%%%%%%%%%%%%%%%%%%%%%%%%%%%%%%%%%%%%%%%%%%%%%%%%%%%%%%%%%%%%%%%%%%%%%%%%%%%%%%%%%%%%%%%%%%%%%%%%%%%%%%%%%%%%%%%%%%%%%%%%%%%

\section{School cohesiveness}\label{Cohesiveness}
In this section, we  want to introduce a scientific definition of  cohesiveness possessed by fish school as a nature of school. We then investigate the relationships between school cohesiveness and obstacle avoiding patterns. These relationships suggest that  avoidance patterns can be used  to visualize the school cohesiveness.

\subsection{Definition and measurement of  school cohesiveness}
Let us first introduce a scientific definition of  cohesiveness for the system \eqref{eq1} of SDEs  in free space. 

\begin{definition}\label{SchoolCohesiveness} {\it School cohesiveness} is the ability of group of fish to form and maintain the $\varepsilon,\theta$-schooling structure against the noise imposed on the school. In other words, how far the group can keep on $\varepsilon,\theta$-schooling as the magnitudes of the noises increase.
\end{definition}

This definition is given in a quantitative form. When $\varepsilon$ and $\theta$ are specified, it is possible to quantitatively measure  the cohesiveness of a group.

Let us next  give some examples of measuring cohesiveness of fish school by numerical methods.

Consider a  group of 50 fish moving in $\mathbb R^2$. The parameters are set as $p=4$, $q=5$, $\alpha=4$, $\beta=1$. The external force functions are taken as  $F_i=-{\mathbf v}_i$ for $i=1,2,\ldots,N$. Initial positions ${\mathbf x}_i(0)$ are randomly located in a suitably small domain with null initial velocities ${\mathbf v}_i(0)={\mathbf 0}$. The magnitude $\sigma_i\equiv \sigma$ is a control parameter of simulation.

 We pick out  20 different trajectories of the Wiener process. For each value $\sigma$, numerical computations for the solution ${\mathbf x}_i(t)$ and ${\mathbf v}_i(t)$ are performed in 20 trials corresponding to these  trajectories. 

Set  $\varepsilon=r=0.5$ and $\theta=0.05$. It is examined whether  or not the states $({\mathbf x}_i(t),{\mathbf v}_i(t)) (i=1,2,\dots,N)$   are in $\varepsilon,\theta$-schooling by fixing $T= 30$. In other words, it is checked that whether or not $n_\varepsilon(t)=1 $ and $\sigma V(t)\leq \theta$ for every $t\geq 30$.

 Starting with sufficiently small $\sigma$, $\sigma$ is then increased with increment step 0.001. When the $\varepsilon,\theta$-schooling structure is broken down at least for one sample trajectory of the Wiener process, the fish are  considered to have lost ability of schooling. The critical value of $\sigma$ which is the largest value of $\sigma$ such that the group is still in  $\varepsilon,\theta$-schooling can  then be found.

\begin{figure}[h!]
\includegraphics[width=1\textwidth]{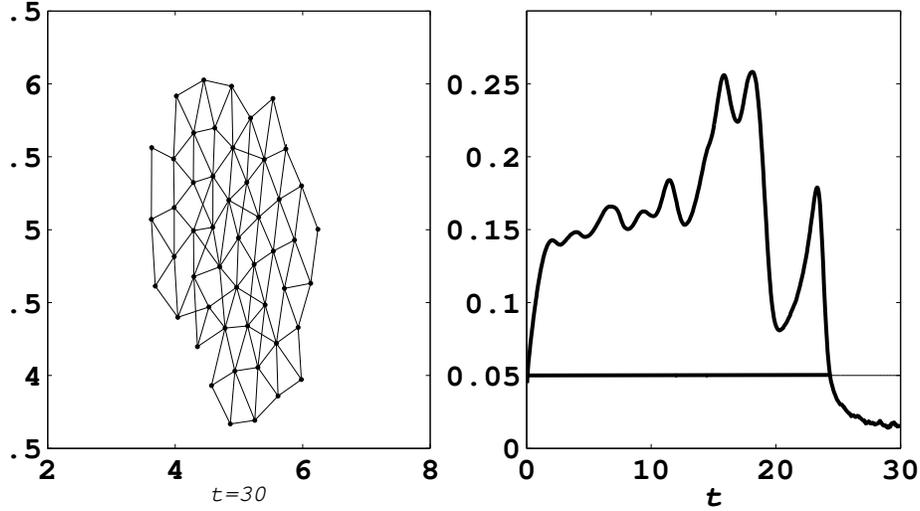}
\caption{$\varepsilon$-graph and $\sigma{\rm VS}(t)$ for $\sigma=0.02$}
\label{fig:2}       
\end{figure}

Fig.~\ref{fig:2} shows that when $\sigma=0.02$,  the group builds up the $\varepsilon,\theta$-schooling structure. Here, the sub-figure on the left shows the $\varepsilon$-graph of the group at time $t=30$, the positions ${\mathbf x}_i(t)$ being drawn by dots, and the edges of graph by lines. The right hand side sub-figure illustrates the variance of velocity as a function of $t$, where the horizontal line represents the level $\theta=0.05$.

\begin{figure}[h!]
	\includegraphics[width=0.75\textwidth]{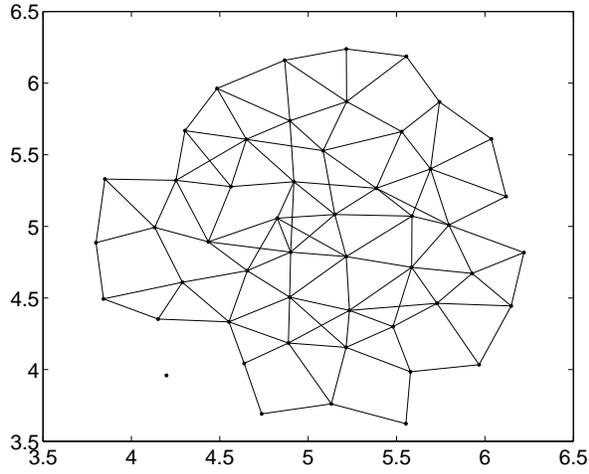}
	\caption{$\varepsilon$-graph for $\sigma=0.06$}
	\label{fig:3}       % Give a unique label
\end{figure}

\begin{figure}[h!]
	\includegraphics[width=0.75\textwidth]{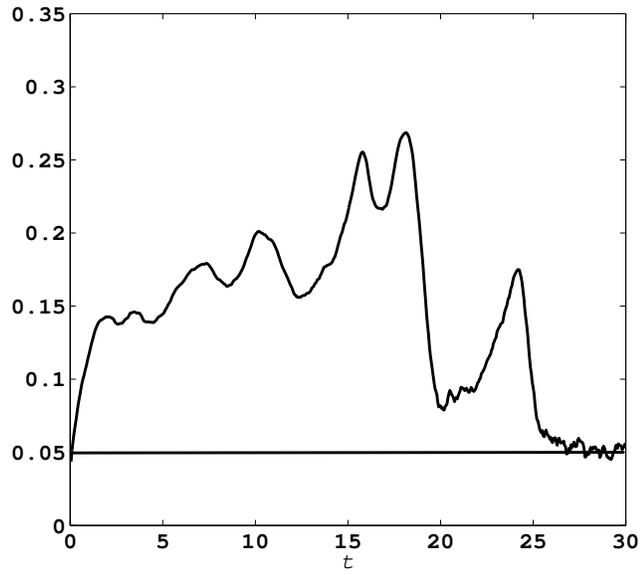}
	\caption{$\sigma {\rm VS}(t)$ for $\sigma=0.069$}
	\label{fig:4}       % Give a unique label
\end{figure}

 The group loses schooling ability when $\sigma=0.06$. In fact,  Fig.~\ref{fig:3} demonstrates that  $n_{\varepsilon}(t)=2$ for $t=30$.  Schooling ability of the group is also lost when  $\sigma=0.069$.  Fig.~\ref{fig:4} highlights that  $\sigma{\rm VS}(t)>\theta$ for $t=30$. 

By these methods,  we can finally compute the critical value $\bar\sigma=0.051$. The group is in $\varepsilon,\theta$-schooling whenever $\sigma \leq \bar \sigma$, whereas it loses schooling ability whenever $\sigma > \bar \sigma$.

\subsection{Relations between exponents $p, q$ and school cohesiveness}  \label{Exponents-Cohesiveness}
Let tune $p$ as $p=4,\, 3.62,\, 3,\, 2$, keeping the relation $q=p+1$ and other parameters as before. Our computations show that the numerical critical value of $\sigma$, namely the cohesiveness, changes as $\bar \sigma=0.051,\, 0.055,\, 0.056, \, 0.063,$ respectively.

The exponent $p$, as explained in Section \ref{MathematicalModels},  shows a degree of range how far the attraction is effective. It is therefore very natural that as $p$ decreases, the attraction range extends and enhances the cohesiveness of group.

\subsection{Relations between critical distance  and school cohesiveness}
Let us tune  $r$ from 0.5, 0.6 and 0.7 by taking $\varepsilon=r$, $p=4$, $q=5$, $\theta=0.05$, $T=30$. The corresponding cohesiveness is  found as $\bar \sigma=0.051, 0.053, 0.054,$ respectively. This means that cohesiveness is enhanced as the critical distance increases.

\subsection{Relationship between school cohesiveness and obstacle avoiding patterns}\label{PatternCohesiveness}
In Subsections \ref{Patterns} and \ref{Exponents-patterns}, we already show that suitable 
tuning of the exponents $p$ and $q=p+1$  provides different obstacle avoiding patterns. It is possible to interpret this fact as follows.

Note that $p$ is concerned with the range of attraction among fish. If $p$ is small, then the attractive force reaches a wide range beyond the critical distance $r$. In contrast, if $p$ is large, then the attraction is only available  in a neighborhood of the disk of radius $r$. 
In Subsection \ref{Exponents-Cohesiveness}, we   verify that when other parameters are fixed, the cohesiveness of school increases as $p$ decreases.

In the numerical examples in Subsection \ref{Patterns}, when $p=2$, the school has very strong cohesiveness and rebounds off the obstacle. When $p=3$, the school  still has  strong cohesiveness and can keep schooling. But the fish are spread on the surface. When $p=3.62$, the school cohesiveness becomes smaller. The fish can no longer keep being in schooling but it is strong enough to reunite the members into a school. When $p=4$, the school cannot keep being in schooling and is separated into two clusters after passing the obstacle.

If these interpretations are reasonable, then the four obstacle avoiding patterns can be used to measure the cohesiveness of school approximately. For example, the cohesiveness can be easily categorized into four classes.

%In conclusion, school cohesiveness can be defined as the ``strength" of group as an internal property to maintain the group against various disturbances. But it is difficult to measure its strength or degree quantitatively. Using the mathematical model, however, we can visualize it by connection with patterns that a school performs in the presence of obstacles. 

%%%%%%%%%%%%%%%%%%%%%%%%%%%%%%%%%%%	Conclusions Section     	%%%%%%%%%%%%%%%%%%%%%%%%%%%%%%%%%%%%%%%%%%%%%%%%%%%%%%%%%%%
\section{Conclusions}\label{Conclusions}
Obstacle avoiding patterns and cohesiveness of fish school have been studied. 
In fact, our mathematical  model for fish schooling in the space with obstacles provided four clear  avoidance patterns:  Rebound, Pullback, Pass and Reunion, and Separation. The shift from one pattern to another due to the change of parameters has also been  investigated.

Furthermore, our definition for  school cohesiveness suggested that the cohesiveness can roughly be measured by using the four types of patterns. In other words, we could visualize the ``strength" of the school by connection with the patterns. 
 Quantitative understanding for the cohesiveness of swarming systems may in turn provide useful information in designing artificial self-organizing systems or intelligent systems \cite{Bonabeau1999,Olfati2003}.

Our attempt may however not be complete, since the emerging patterns depend not only on the internal parameters  but also on the environmental conditions, for instance the radius of obstacle.   If the geometrical diameter of school is relatively larger than the obstacle's radius, then the school striking the obstacle may pass over and reunify. To the contrary, if the school's diameter is relatively smaller than the radius, then the school may rebound off. 

 %Our results may suggest possibility to study the school cohesiveness using mathematical models.

\section*{Acknowledgements}
The authors heartily express their gratitude to the referees of this paper for making useful and 
 constructive comments on the style of paper. 
The work of the second author was supported by JSPS KAKENHI Grant Number 20140047.

%\section*{References}

\end{document}